\documentclass[aos,citesort,dvips]{arximspdf}
\usepackage{graphics}
\doi{10.1214/10-AOS791}
\volume{38}
\issue{4}
\pubyear{2010}
\firstpage{2559}
\lastpage{2586}

\makeatletter

\newtheorem{theorem}{Theorem}[section]
\newtheorem{lemma}{Lemma}[section]

\makeatother

\begin{document}
\begin{frontmatter}

\title{Fractals with point impact in functional linear regression}
\runtitle{Fractals with point impact}

\begin{aug}
\author[A]{\fnms{Ian W.} \snm{McKeague}\thanksref{t1}\ead[label=e1]{im2131@columbia.edu}\corref{}} and
\author[B]{\fnms{Bodhisattva} \snm{Sen}\thanksref{t2}\ead[label=e2]{bs2528@columbia.edu}}
\runauthor{I. W. McKeague and B. Sen}
\affiliation{Columbia University}
\address[A]{Department of Biostatistics\\
Columbia University\\
722 West 168th Street, 6th Floor\\
New York, New York 10032\\
USA\\
\printead{e1}}
\address[B]{Department of Statistics\\
Columbia University\\
1255 Amsterdam Avenue, 10th Floor\\
New York, New York 10027\\
USA\\
\printead{e2}}
\end{aug}

\thankstext{t1}{Supported by NSF Grant DMS-08-06088.}
\thankstext{t2}{Supported by NSF Grant DMS-09-06597.}

\received{\smonth{7} \syear{2009}}
\revised{\smonth{12} \syear{2009}}

%
\begin{abstract}
This paper develops a point impact linear regression model in which
the trajectory of a continuous stochastic process, when evaluated at a
sensitive time point, is associated with a scalar response. The
proposed model complements and is more interpretable than the
functional linear regression approach that has become popular in recent
years. The trajectories are assumed to have fractal (self-similar)
properties in common with a fractional Brownian motion with an unknown
Hurst exponent. Bootstrap confidence intervals based on the
least-squares estimator of the sensitive time point are developed.
Misspecification of the point impact model by a functional linear model
is also investigated. Non-Gaussian limit distributions and rates of
convergence determined by the Hurst exponent play an important role.
\end{abstract}

%

\end{frontmatter}

\section{Introduction}\label{sec1}

This paper investigates a linear regression model involving a scalar
response $Y$ and a predictor given by the value of the trajectory of a
continuous stochastic process $X=\{X(t)$, $t\in[0,1]\}$ at some
unknown time point. Specifically, we consider the \textit{point impact}
linear regression model
%
%
\begin{equation}\label{workmodel}
Y= \alpha+\beta X(\theta) + \varepsilon
\end{equation}
and focus on the time point $\theta\in(0,1)$ as the target parameter
of interest. The intercept $\alpha$ and the slope $\beta$ are
scalars, and the error $\varepsilon$ is taken to be independent of $X$,
having zero mean and finite variance $\sigma^2$. The complete
trajectory of $X$ is assumed to be observed (at least on a fine enough
grid that it makes no difference in terms of accuracy), even though the
model itself only involves the value of $X$ at $\theta$, which
represents a ``sensitive'' time point in terms of the relationship to
the response. The main aim of the paper
is to show that the precision of estimation of $\theta$ is driven by
\textit{fractal} behavior in $X$,
and to develop valid inferential procedures that adapt to a broad range
of such behavior.
Our model could easily be extended in various ways, for example, to
allow multiple sensitive time points or further covariates, but, for
simplicity, we restrict attention to (\ref{workmodel}).

\begin{figure}

\includegraphics{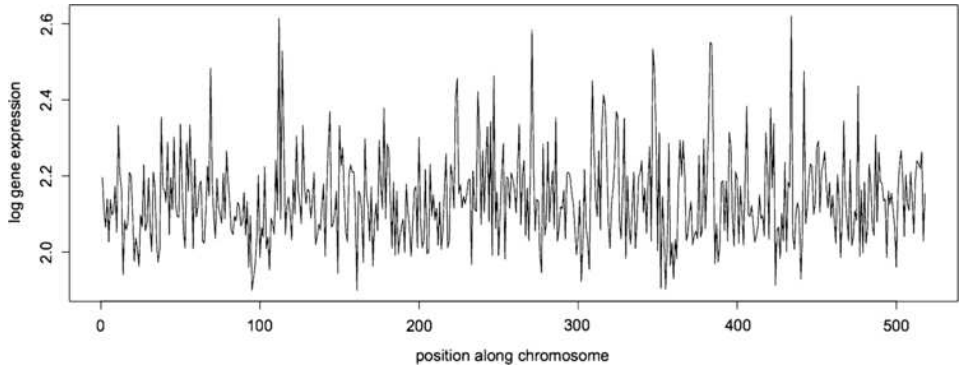}

\caption{Log gene expression at 518 loci along chromosome 17 in tissue
from a breast cancer patient.}\label{fig1}
\end{figure}

Our motivation for developing this type of model arises from
genome-wide expression studies that measure the activity of numerous
genes simultaneously.
In these studies, it is of interest to locate genes showing activity
that is associated with clinical outcomes. Emilsson et al.
\cite{emilsson}, for example, studied gene expression levels at over 24,000
loci in samples of adipose tissue to identify genes correlated with
body mass index and other obesity-related outcomes. Gruvberger-Saal et
al. \cite{gruvberger} used gene expression profiles from the tumors of
breast cancer patients to predict estrogen receptor protein
concentration, an important prognostic marker for breast tumors; see
also \cite{buness}. In such studies, the gene expression profile
across a chromosome can be regarded a functional predictor, and a gene
associated with the clinical outcome is identified by its base pair
position $\theta$ along the chromosome; see Figure \ref{fig1}. Our
aim here is to develop a method of estimating a confidence interval for
$\theta$, leading to the identification of chromosomal regions that
are potentially useful for diagnosis and therapy. Although there is
extensive statistical literature on gene expression data, it is almost
exclusively concerned with multiple testing procedures for detecting
differentially expressed genes; see, for example, \cite{dula,salas}.

%

Gene expression profiles (as in Figure \ref{fig1}) clearly display
fractal behavior, that is,
self-similarity over a range of scales. Indeed, fractals often arise when
spatiotemporal patterns at higher levels emerge from localized
interactions and selection processes acting at lower levels, as with
gene expression activity.
Moreover, the recent discovery \cite{l-a} that chromosomes are folded as
``fractal globules,'' which can easily unfold during
gene activation, also helps explain the fractal appearance of gene
expression profiles.

A basic stochastic model for fractal phenomena is provided by
fractional Brownian motion (fBm) (see \cite{manness}), in which the
so-called Hurst exponent $H\in[0,1]$
calibrates the scaling of the self-similarity
and provides a natural measure of trajectory roughness.
It featured prominently in the pioneering work of Beno\^{i}t
Mandelbrot, who
stated (\cite{man}, page 256) that fBm provides ``the most manageable
mathematical environment I can think of (for representing fractals).''
For background on fBm from a statistical modeling point of view, see
\cite{gn}.

The key issue to be considered in this paper is how to construct a
confidence interval for the true sensitive time point $\theta_0$
based on its least squares estimator $\hat\theta_n$, obtained by
fitting model (\ref{workmodel}) from a sample of size $n$,
%
%
\begin{equation}\label{lsq}
(\hat\alpha_n, \hat\beta_n,\hat\theta_n)
=\mathop{\arg\min}_{\alpha, \beta,\theta} \sum_{i=1}^n[Y_i- \alpha-\beta
X_i(\theta)]^2.
\end{equation}
We show that, when $X$ is fBm, both the rate of convergence $r_n$ and
limiting distribution
of $\hat\theta_n$ depend on $H$. In addition, we construct bootstrap
confidence intervals for $\theta_0$ that \textit{do not require knowledge
of} $H$. This facilitates applications (e.g., to gene expression data)
in which the type of fractal behavior is not known in advance; the
trajectory in Figure \ref{fig1} has an estimated Hurst exponent of
about 0.1, but it would be very difficult to estimate precisely using
data in a small neighborhood of $\hat\theta_n$, so a bootstrap
approach becomes crucial.
We emphasize that nothing about the distribution of $X$ is used in the
construction of the estimators or the bootstrap confidence intervals;
the fBm assumption will only be utilized to study the large sample
properties of these procedures.
Moreover, our main results will make essential use of the fBm
assumption only \textit{locally}, that is, in a small neighborhood
of~$\theta_0$.

The point impact model (\ref{workmodel}) can be regarded as a simple
working model that provides interpretable information about the
influence of $X$ at a specific location (e.g., a genetic locus). Such
information cannot be extracted using the standard functional linear
regression model \cite{rs} given by
%
%
\begin{equation}\label{funlinear}
Y=\alpha+ \int_0^1 f(t) X(t) \,dt + \varepsilon,
\end{equation}
where $f$ is a continuous function and $\alpha$ is an intercept,
because the influence of $X(t)$ is spread continuously across $ [0,1]$
and point-impact effects are excluded. In the gene expression context,
if only a few genes are predictive of $Y$, then a model of the form
(\ref{workmodel}) would be more suitable than (\ref{funlinear}),
which does not allow $f$ to have infinite spikes. In general, however,
a continuum of locations is likely to be involved (as well as point-impacts),
so it is of interest to study the behavior of $\hat\theta_n$ in
misspecified settings in which
the data arise from combinations of (\ref{workmodel}) and (\ref{funlinear}).

Asymptotic results for the least squares estimator (\ref{lsq}) in the
correctly specified
setting are presented in Section \ref{LSE}. In Section \ref{BootsCI}
it is shown that the residual bootstrap is consistent for
the distribution of $\hat\theta_n$, leading to the construction of
valid bootstrap confidence intervals without knowing $H$. The
nonparametric bootstrap is shown to be inconsistent
in the same setting. The effect of misspecification is discussed in
Section \ref{MissFLM}. A two-sample problem
version of the point impact model is discussed in Section \ref{two}.
Some numerical examples are presented in Section \ref{NumEx}, where we
compare the proposed bootstrap confidence interval with Wald-type
confidence intervals (in which $H$ is assumed to be known); an
application to gene expression data is also discussed. Concluding
remarks appear in Section \ref{Discuss}. Proofs are placed in Section
\ref{proofs}.

\section{Least squares estimation of the sensitive time point}\label{LSE}
Throughout we take $X$ to be a fBm with Hurst exponent $H$, which, as
discussed earlier, controls the roughness of the trajectories.
We shall see in this section that the rate of convergence of $\hat
\theta_n$ can be expressed explicitly in terms of $H$.

First we recall some basic properties of fBm. A (standard) fBm with
Hurst exponent $H\in(0,1]$ is a Gaussian process $B_H=\{ B_H(t)$,
$t\in\mathbb{R}\}$ having continuous sample paths, mean zero and
covariance function
%
%
\begin{equation}\label{covar}
\operatorname{Cov}\{ B_H(t),B_H(s)\} = \tfrac{1}{2} ( |t|^{2H}
+|s|^{2H}-|t-s|^{2H} ).
\end{equation}
By comparing their mean and covariance functions, $B_H(a t) \stackrel
{d}{=} a ^{H}B_H(t)$ as processes, for all $a>0$ (self-similarity).
Clearly, $B_{1/2}$ is a two-sided Brownian motion, and $B_1$ is a
random straight line: $B_1(t)=tZ$ where $Z\sim N(0,1)$. The increments
are negatively correlated if $H<1/2$, and positively correlated if
$H>1/2$. Increasing $H$ results in smoother sample paths.

Suppose $(X_i,Y_i), i=1,\ldots,n$, are i.i.d. copies of $(X,Y)$
satisfying the model (\ref{workmodel}).
The unknown parameter is $\eta=(\alpha,\beta,\theta)\in\Xi=
\mathbb{R}^2 \times[0,1]$, and its true value is denoted $\eta_0 =
(\alpha_0,\beta_0, \theta_0)$. The following conditions are needed:

\begin{itemize}[(A3)]
\item[(A1)] $X$ is a fBm with Hurst exponent $H \in(0,1)$.

\item[(A2)]$0<\theta_0<1$ and $\beta_0 \ne0$.

\item[(A3)]$E |\varepsilon|^{2 + \delta} < \infty$ for some $\delta
> 0$.
\end{itemize}

The construction of the least squares estimator $\hat\eta_n=(\hat
\alpha_n, \hat\beta_n,\hat\theta_n)$, defined by (\ref{lsq}),
does not involve any assumptions about the distribution of the
trajectories, whereas the asymptotic behavior does. Our first result
gives the consistency and asymptotic distribution of $\hat\eta_n$
under the above assumptions.
\begin{theorem} \label{fbmcase} If \textup{(A1)} and \textup{(A2)} hold, then
$\hat\eta_n$ is consistent, that is, $\hat\eta_n \stackrel
{P}{\rightarrow} \eta_0$.
If \textup{(A3)} also holds, then
%
%
\begin{eqnarray}\label{eq:LimDist}
\zeta_n & \equiv& \bigl( \sqrt{n}(\hat\alpha_n - \alpha_0), \sqrt
{n}(\hat\beta_n - \beta_0), n^{{1/(2H)}} (\hat\theta_n-\theta_0)
\bigr) \nonumber\\[-8pt]\\[-8pt]
& \stackrel{d}{\rightarrow} & \biggl( \sigma Z_1, |\theta_0|^{-H} \sigma
Z_2, \mathop{\arg\min}_{t \in\mathbb{R}} \biggl\{ 2 \frac{\sigma}{|\beta_0|}
B_H(t) + |t|^{2H} \biggr\} \biggr) \equiv\zeta,\nonumber
\end{eqnarray}
where $Z_1$ and $Z_2$ are i.i.d. $N(0,1)$, independent of the fBm $B_H$.
\end{theorem}

\subsection*{Remarks}
\begin{enumerate}
\item It may come as a surprise that the convergence rate of $\hat
\theta_n$ increases as $H$ decreases, and becomes arbitrarily fast as
$H\to0$. A heuristic explanation is that fBm ``travels further'' with
a smaller $H$, so independent trajectories of $X$ are likely to ``cover
different ground,'' making it easier to estimate $\theta_0$. In a
nutshell, the smaller the Hurst exponent, the better the design.

\item It follows from (a sight extension of) Lemmas 2.5 and 2.6 of Kim
and Pollard \cite{kp} that the third component of $\zeta$ is well defined.

\item Using the self-similarity of fBm, the asymptotic distribution of
$\hat\theta_n$ can be expressed as the distribution of
%
%
\begin{equation}
\label{delta}
\Delta\equiv\biggl( \frac{\sigma}{|\beta_0|} \biggr)^{1/H} \mathop{\arg\min}_{t \in
\mathbb{R}} \bigl( B_H(t) +|t|^{2H}/2 \bigr).
\end{equation}
This distribution does not appear to have been studied in the
literature except for $H=1/2$ and $H=1$ (standard normal). When
$H=1/2$, $X$ is a standard Brownian motion and the limiting
distribution is given in terms of a two-sided Brownian motion with a
triangular drift. Bhattacharya and Brockwell \cite{bb} showed that
this distribution has a density that can be expressed in terms of the
standard normal distribution function. It arises frequently in
change-point problems under contiguous asymptotics \cite{y,s,ms}.

\item From the proof, it can be seen that the essential assumptions on
$X$ are the self-similarity and stationary increments properties in
some neighborhood of~$\theta_0$, along with the trajectories of $X$
being Lipschitz of all orders less than $H$. Note that any Gaussian
self-similar process with stationary increments and zero mean is a fBm
(see, e.g., Theorem 1.3.3 of \cite{em}).

\item The trajectories of fBm are nondifferentiable when $H<1$, so the
usual technique of Taylor expanding the criterion function about
$\theta_0$ does not work and a more sophisticated approach is required
to prove the result.

\item Note that $(\hat\alpha_n,\hat\beta_n)$ has the same limiting
behavior as though $\theta_0$ is known, and $\hat\theta_n$ and
$(\hat\alpha_n,\hat\beta_n)$ are asymptotically independent.

\item
The result is readily extended to allow for additional covariates [cf.
(\ref{blinear})], which is often
important in applications. The limiting distribution of $\hat\theta
_n$ remains the same, and the other regression
coefficient estimates have the same limiting behavior as though $\theta
_0$ is known.

\item Note that the assumption $\beta_0 \ne0$ is crucial for the
theorem to hold.
When \mbox{$\beta_0=0$}, the fBm does not influence the response at all and
$\hat\theta_n$ contains no information about $\theta_0$.

\end{enumerate}

\section{Bootstrap confidence intervals}\label{BootsCI}
In general, a valid Wald-type confidence interval for $\theta_0$ would
at least need a consistent estimator of the Hurst exponent $H$, which
is a nuisance parameter in this problem. Unfortunately, however,
accurate estimation of $H$ is difficult and quite often unstable.
Bootstrap methods have been widely applied to avoid issues of nuisance
parameter estimation, and they work well in problems with $\sqrt
{n}$-rates; see \cite{bickel81,singh81,shao95} and the references
therein. In this section we study the consistency properties of two
bootstrap methods that arise naturally in our setting. One of these
methods leads to a valid confidence interval for $\theta_0$ without
requiring any knowledge of $H$.

\subsection{Preliminaries}\label{IntroBoots}
We start\vspace*{-1pt} with a brief review of the bootstrap. Given a sample ${\mathbf
Z}_n=\{Z_1, Z_2, \ldots, Z_n\}\stackrel{\mathrm{i.i.d.}}{\sim} L$ from an
unknown\vspace*{1pt} distribution $L$,
suppose that the distribution function, $F_n$, say, of some random
variable $R_n \equiv R_n(\mathbf{Z}_n, L)$, is of interest; $R_n$ is
usually called a root and it can in general be any measurable function
of the data and the distribution $L$. The bootstrap method can be
broken into three simple steps:
\begin{longlist}
\item Construct an estimator $\hat{L}_n$ of $L$ from ${\mathbf Z}_n$.

\item Generate ${\mathbf Z}_n^{*} =\{Z_1^{*},\ldots,
Z_{n}^{*}\} \stackrel{\mathrm{i.i.d.}} {\sim} \hat{L}_n$ given
${\mathbf Z}_n$.

\item Estimate $F_n$ by $F_{n}^*$, the conditional c.d.f. of
$R_n({\mathbf Z}_n^{*},\hat{L}_n)$ given ${\mathbf Z}_n$. 
\end{longlist}
Let $d$ denote the L\'{e}vy metric or any other metric metrizing weak
convergence of distribution functions. We say that $F_{n}^*$ is weakly
consistent if $d(F_n,F_n^*)\stackrel{P}{\rightarrow} 0$;
if $F_{n}$ has a weak limit $F$, this is equivalent to $F_{n}^*$
converging weakly to $F$ in probability.

The choice of $\hat{L}_n$ mostly considered in the literature is the
empirical distribution. Intuitively, an $\hat L_n$ that mimics the
essential properties (e.g., smoothness) of the underlying distribution
$L$ can be expected to perform well. In most situations, the empirical
distribution of the data is a good estimator of $L$, but in some
nonstandard situations it may fail to capture some of the important
aspects of the problem, and the corresponding bootstrap method can be
suspect. The following discussion illustrates this phenomenon (the
inconsistency when bootstrapping from the empirical distribution of the
data) when $\Delta_n \equiv n^{1/(2H)} (\hat\theta_n - \theta_0)$
is the random variable of interest.

\subsection{Inconsistency of bootstrapping pairs}\label{sec32}
In a regression setup there are two natural ways of bootstrapping:
bootstrapping pairs (i.e., the nonparametric bootstrap) and
bootstrapping residuals (while keeping the predictors fixed). We show
that bootstrapping pairs (drawing $n$ data points with replacement from
the original data set) is inconsistent for $\theta_0$.
\begin{theorem}\label{thm:InConsBoots}
Under conditions \textup{(A1)--(A3)}, the nonparametric bootstrap is
inconsistent for estimating the distribution of $\Delta_n$, that is,
$\Delta_n^* \equiv n^{1/(2H)} (\hat\theta_n^* - \hat\theta_n)$,
conditional on the data, does not converge in distribution to $\Delta$
in probability, where $\Delta$ is defined by (\ref{delta}).
\end{theorem}

\subsection{Consistency of bootstrapping residuals}\label{sec33}
Another bootstrap procedure is to use the form of the assumed model
more explicitly to draw the bootstrap samples: condition on the
predictor $X_i$ and generate its response as
%
%
\begin{equation}\label{eq:DefYi*}
Y_i^* = \hat\alpha_n + \hat\beta_n X_i(\hat\theta_n) +
\varepsilon_i^*,
\end{equation}
where the $\varepsilon_i^*$ are conditionally i.i.d. under the empirical
distribution of the centered
residuals $\hat\varepsilon_i - \bar\varepsilon_n$, with $\hat
\varepsilon_i
= Y_i - \hat\alpha_n - \hat\beta_n X_i(\hat\theta_n)$ and $\bar
\varepsilon_n = \sum_{i=1}^n \hat\varepsilon_i/n$. Let $\hat\alpha_n^*,
\hat\beta_n^*$ and $\hat\theta_n^*$ be the estimates of the unknown
parameters obtained from the bootstrap sample. We approximate the
distribution of $\zeta_n$ [see (\ref{eq:LimDist})] by the conditional
distribution of
\[
\zeta_n^* \equiv\bigl[ \sqrt{n}(\hat\alpha_n^* - \hat\alpha_n),
\sqrt{n}(\hat\beta_n^* - \hat\beta_n), n^{1/(2H)}(\hat\theta_n^*
- \hat\theta_n) \bigr],
\]
given the data.
\begin{theorem}\label{ConsBoots}
Under conditions \textup{(A1)--(A3)}, the above procedure of bootstrapping
residuals is consistent for estimating the distribution of $\zeta_n$,
that is, $\zeta_n^* \stackrel{d}{\rightarrow} \zeta$, in
probability, conditional on the data.
\end{theorem}

We now use the above theorem to construct a valid confidence interval
(CI) for $\theta_0$ that does not require any knowledge of $H$. Let
$q^*_{\alpha}$ be the $\alpha$-quantile of the conditional
distribution of $(\hat\theta_n^* - \hat\theta_n)$ given the data,
which can be readily obtained via simulation and does not involve the
knowledge of any distributional properties of $X$. The proposed
approximate $(1 - 2 \alpha)$-level bootstrap CI for $\theta_0$ is
then given by
\[
\mathcal{C}_n = [\hat\theta_n - q^*_{1 - \alpha}, \hat\theta_n -
q^*_{\alpha}].
\]
Under the assumptions of Theorem \ref{ConsBoots}, the coverage
probability of this CI is
\begin{eqnarray*}
P \{\theta_0 \in\mathcal{C}_n \}
& = & P\bigl\{n^{1/(2H)} q^*_{\alpha} \le\Delta_n \le n^{1/(2H)} q^*_{1 -
\alpha}\bigr\} \\
& \approx& P^*\bigl\{n^{1/(2H)} q^*_{\alpha} \le\Delta_n^* \le
n^{1/(2H)} q^*_{1 - \alpha}\bigr\} \\
& = & P^*\{q^*_{\alpha} \le\hat\theta_n^* - \hat\theta_n \le
q^*_{1 - \alpha}\} \\
&=& 1 - 2 \alpha,
\end{eqnarray*}
where $P^*$ denotes the bootstrap distribution given the data, and we
have used the fact that the supremum distance between the relevant
distribution functions of $\Delta_n$ and $\Delta_n^*$ is
asymptotically negligible. The key point of this argument is that
$\Delta_n$ and $\Delta_n^*$ have the same normalization factor
$n^{1/(2H)}$ and, thus, it ``cancels'' out.
CIs for $\alpha_0$ and $\beta_0$ can be constructed in a similar fashion.

\subsection{Discussion}\label{sec34}
In nonparametric regression settings, dichotomies in the behavior of
different bootstrap methods are well known, for example, when using the
bootstrap to calibrate omnibus goodness-of-fit tests for parametric
regression models; see \cite{hm,vkgms,neu} and references therein. A
dichotomy in the behavior of the two bootstrap methods, however, is
surprising in a linear regression model. This illustrates that in
problems with nonstandard asymptotics, the usual nonparametric
bootstrap might fail, whereas a resampling procedure that uses some
particular structure of the model can perform well. The improved
performance of bootstrapping residuals will be confirmed by our
simulation results in Section \ref{NumEx}.

The difference in the behavior of the two bootstrap methods can be
explained as follows. As in any M-estimation problem, the standard
approach is to study the criterion (objective) function being
optimized, in a neighborhood of the target parameter, by splitting it
into an empirical process and a drift term. The drift term has
different behavior for the two bootstrap methods: while bootstrapping
pairs, it does not converge, whereas the bootstrapped residuals are
conditionally independent of the predictors, and the drift term
converges. This highlights the importance of designing the bootstrap to
accurately replicate the structure in the assumed model. A~more
technical explanation is provided in a remark following the proof of
Theorem~\ref{ConsBoots}.

Other types of resampling (e.g., the $m$-out-of-$n$ bootstrap, or
subsampling) could be applicable,
but such methods require knowledge of the rate of convergence, which
depends on the unknown $H$.
Also, these methods require the choice of a tuning parameter, which is
problematic in practice.
However, the residual bootstrap is consistent, easy to implement, and
does not require the
knowledge of $H$ and the estimation of a tuning parameter.

The inconsistency of the nonparametric bootstrap casts some doubt on
its validity for checking the stability of variable selection results
in high-dimensional regression problems (as is common practice).
Indeed, it suggests that more care (in terms of more explicit use of
the model) is needed in the choice of a bootstrap method in such settings.

\section{Misspecification by a functional linear model}\label{MissFLM}
The point impact model
cannot capture effects that are spread out over the domain of the
trajectory, for example, gene expression
profiles for which the effect on a clinical outcome involves complex
interactions between numerous genes.
Such effects, however, may be represented by a functional linear model,
and we now examine how the limiting behavior of $\hat\theta_n$
changes when the data arise in this way.


\subsection{Complete misspecification}\label{sec41}
In this case we treat (\ref{workmodel}) as the working model (for fitting the data),
but view this model as being completely misspecified in the sense that
the data are generated from the functional linear model (\ref{funlinear}).
For simplicity, we set $\alpha=0$ and $\beta=1$ in the working model,
and set
$\alpha= 0$ in the true functional linear model. The least squares
estimator $\hat\theta_n$
now estimates the minimizer $\theta_0$ of
\[
\mathbb{M}(\theta) \equiv E[Y-X(\theta)]^2=\sigma^2 +E \biggl[ \int_0^1
f(t)X(t) \,dt - X(\theta) \biggr]^2
\]
and the following result gives its asymptotic distribution.
%
%
\begin{theorem}\label{miss}
Suppose that \textup{(A1)} and \textup{(A3)} hold, and that $\mathbb{M}(\theta)$ has a unique
minimizer and is twice-differentiable at $0<\theta_0<1$. Then, in the
misspecified case,
\[
n^{{1/(4-2H)}}(\hat\theta_n-\theta_0) \stackrel{d}{\rightarrow}
\mathop{\arg\min}_{t\in\mathbb{R}}\bigl(2aB_H(t) +bt^2\bigr),
\]
where $a^2= \mathbb{M}(\theta_0)$ and $ b =
\mathbb{M}''(\theta_0)/2$.
\end{theorem}

\subsection*{Remarks}
\begin{enumerate}
\item Here the rate of convergence reverses itself from the correctly
specified case: the convergence rate now decreases as $H$ decreases,
going from a parametric rate of $n^{1/2}$ when $H\to1$, to as slow as
$n^{1/4}$ as $H\to0$. A heuristic explanation
is that roughness in $X$ now amounts to measurement error (which
results in a slower rate) as the
fluctuations of $X$ are smoothed out in the true model.

\item
In the case of Brownian motion trajectories ($H=1/2$), note that
$\mathbb{M}(\theta)=
\theta- 2\int_0^1f(t)\min(t,\theta) \,dt+\mathrm{const}$, the normal
equation is
%
%
\begin{equation}
\label{normal}
\mathbb{M}'(\theta)=1 -2\int_\theta^1f(t) \,dt=0
\end{equation}
%
and $\mathbb{M}''(\theta)=2f(\theta)$.

\item Also in the case $H=1/2$, the limiting distribution is given in
terms of two-sided Brownian motion with a parabolic drift, and was
investigated originally by Chernoff \cite{c} in connection with the
estimation of the mode of a distribution, and shown to have a density
(as the solution of a heat equation). The Chernoff distribution arises
frequently in monotone function estimation settings; 
Groeneboom and Wellner \cite{gw} introduced various algorithms for
computation of its distribution function and quantiles.
\end{enumerate}

\subsection{Partial misspecification}\label{sec42}
The nonparametric functional
linear model (\ref{funlinear}) can be combined with (\ref{workmodel})
to give the semiparametric model
%
%
\begin{equation}\label{partlinear}
Y=\alpha+ \beta X(\theta) + \int_0^1 f(t) X(t) \,dt + \varepsilon,
\end{equation}
which allows $X$ to have both a point impact and an influence that is
spread out continuously in time. When $f=0$, this model reduces to the
point impact model; when $\beta=0$, to the functional linear model. In
this section we examine the behavior of $\hat\theta_n$ when the
working model is (\ref{workmodel}), as before, but the data are now
generated from~(\ref{partlinear}).

For simplicity, suppose that $\alpha=0$ and $\beta=1$ in both the
working point impact model and in the true
model (\ref{partlinear}). Denote the true value of $\theta$ by
$\theta_0\in(0,1)$. It can then be shown that $\hat\theta_n$ is robust
to small levels of misspecification, that is, it consistently estimates
$\theta_0$ with the same rate of convergence as in the correctly
specified case. Indeed, $\hat\theta_n $
targets the minimizer of the criterion function
\[
\mathbb{M}(\theta) = E[Y - X(\theta) ]^2=|\theta- \theta_0|^{2H}
-\int
_0^1f(t)[t^{2H} + \theta^{2H} - |\theta-t|^{2H}] \,dt + \mathrm{const}.
\]
Provided $\int|f|$ is sufficiently small, the derivative of $\mathbb{M}$
will be negative over the interval $(0,\theta_0)$ and positive over
$(\theta_0,1)$, so $\mathbb{M}$ is minimized at $\theta_0$.
It is then possible to extend Theorem \ref{fbmcase} to give
%
%
\begin{equation}
\label{asympmiss}
n^{{1/(2H)}}(\hat\theta_n-\theta_0) \stackrel{d}{\rightarrow}
a^{1/H} \mathop{\arg\min}_{t\in\mathbb{R}}\bigl(B_H(t) +|t|^{2H}/2\bigr),
\end{equation}
where $a\ge\sigma$ is defined in the statement of Theorem \ref{miss}.
This shows that the effect of partial misspecification is a simple
inflation of the variance [cf. (\ref{delta})], without any change in
the form of the limit distribution.

It is also of interest to estimate $\theta_0$ in a way that adapts to
any function
$f$ (i.e., sufficiently smooth) in this semiparametric setting. This
could be done, for example, by approximating $f$ by a finite B-spline basis
expansion of the form $f_m(t)=\sum_{j=1}^m\beta_j\phi_j(t)$, and
fitting the working model
%
%
\begin{equation}\label{blinear}
Y=\alpha+ \beta X(\theta) + \sum_{j=1}^m\beta_jZ_j + \varepsilon,
\end{equation}
where $Z_j=\int_0^1 \phi_j(t) X(t) \,dt$ are additional covariates with
regression coefficients~$\beta_j$; the resulting least squares
estimator $\tilde\theta_n$ can then be used
as an estimate of $\theta_0$ of $\theta$.
For the working model (\ref{blinear}), the misspecification is
$f-f_m$, which will be small if $m$ is sufficiently large.
Therefore, based on our previous discussion,
$\tilde\theta_n$~will satisfy a result of the form (\ref{asympmiss}); in particular,
$\tilde\theta_n$ will exhibit the fast
$n^{1/(2H)}$-rate of convergence. Note that for this result to hold,
$m$ can be fixed and does not need to tend to infinity with the sample size.

\section{Two-sample problem}\label{two}

In this section we discuss a variation of the point impact regression
model in which the response takes just
two values (say $\pm1$). This is of interest, for example, in
case-control studies in which gene-expression data are available for a
sample of cancer patients
and a sample of healthy controls, and the target parameter is the locus of
a differentially expressed gene.

Suppose we have two independent samples of trajectories $X$, with $n_1$
trajectories from class 1, and $n_2$ trajectories from class $-1$, for
a total sample size of $n=n_1+n_2$. It is assumed that $\rho
=n_1/n_2>0$ remains fixed, and the $j$th sample satisfies the model
\[
X_{ij}(t)=\mu_j(t) +\varepsilon_{ij}(t),\qquad j=1,2,
\]
where $\varepsilon_{ij}$, $i=1,\ldots,n_j$ are i.i.d. fBms with Hurst
exponent $H\in(0,1)$, and $\mu_j(t)$ is an unknown mean function,
assumed to be continuous. The treatment effect $\mathbb{M}(t)=\mu
_1(t)-\mu
_2(t)$ is taken to have a point impact in the sense of having a unique
maximum at $\theta_0\in(0,1)$;
minima can of course be treated in a similar fashion.
The least squares estimator of the sensitive time point now becomes
%
%
\begin{equation}
\label{lsq2}
\hat\theta_n
=\mathop{\arg\max}_{\theta} \{ \bar X_{1}(\theta) - \bar X_{2}(\theta)
\},
\end{equation}
where $ \bar X_{j}(\theta) =\sum_{i=1}^{n_j}X_{ij} (\theta)/{n_j}$
is the sample mean for class $j$. Although a studentized version
(normalizing the the difference of the sample means by
a standard error estimate) might be preferable in some cases, with
small or unbalanced samples, say, to keep the discussion
simple, we restrict attention to $\hat\theta_n$.
The empirical criterion function $\mathbb{M}_n(\theta)= \bar
X_{1}(\theta)
- \bar X_{2}(\theta)$
converges
uniformly to $\mathbb{M}(\theta)$ a.s. (by the Glivenko--Cantelli
theorem), so $\hat\theta_n$ is a consistent estimator of $\theta_0$.

As before, our objective is to find a confidence interval for $\theta
_0$ based on $\hat\theta_n$ under appropriate
conditions on the treatment effect.
Toward this end,
we need an assumption on the degree of smoothness of the
treatment effect at $\theta_0$
in terms of an exponent $0<S\le1$:
\[
\mathbb{M}(\theta)=\mathbb{M}(\theta_0)-c|\theta-\theta_0|^{2S}
+o(|\theta
-\theta_0|^{2S})
\]
as $ \theta\to\theta_0$, where $c>0$. If $\mathbb{M}$ is
twice-differentiable at $\theta_0$, then
this assumption holds only with $S=1$; for it to hold for some $S<1$, a
cusp is needed.
When the smoothness of the treatment effect and the fBm match, that is, $S=H$,
the rate of convergence of $\hat\theta_n$ is $n^{1/(2H)}$, as before, and
$\hat\theta_n$ has a nondegenerate limit distribution of the same
form as in Theorem \ref{fbmcase}:
%
%
\begin{equation}
\label{2samp}
n_1^{1/(2H)} (\hat\theta_n-\theta_0) \stackrel{d}{\rightarrow}
\mathop{\arg\min}_{t\in\mathbb{R}}\bigl\{ \bigl(1+\sqrt\rho\bigr) B_H(t) +c|t|^{2H}\bigr\}.
\end{equation}
The key step in the proof (which is simpler than in the regression
case) is given at the end of Section \ref{proofs}.


\section{Numerical examples}\label{NumEx}
In this section we report some numerical results based on trajectories
from fBm simulations and from gene expression data.

We first consider a correctly specified example as in Section \ref
{LSE} and study the behavior of CIs for the sensitive time-point
$\theta_0$ using
the two bootstrap based methods, and compare them with the
$100(1-\alpha)\%$
Wald-type CI
%
%
\begin{equation}\label{eq:ci}
\hat\theta_n\pm\biggl(\frac{\hat\sigma_n}{|\hat\beta_n| \sqrt n}
\biggr)^{1/H}z_{H,\alpha/2}
\end{equation}
with $H$ assumed to be known.
Here $\hat\sigma_n$ is the sample standard deviation of the
residuals, and $z_{H,\alpha}$ is the upper $\alpha$-quantile of $\arg
\min_{t\in\mathbb{R}}( B_H(t) +|t|^{2H}/2)$. In practice,
$H$ needs to be estimated to apply (\ref{eq:ci}). Numerous estimators
of $H$ based on a single realization of $X$ have been proposed in the
literature \cite{b,coeur}, although observation at fine time scales is
required for such estimators to work well, and it is not clear that
direct plug-in would be satisfactory.
The quantiles $z_{H,\alpha/2}$ needed to compute the Wald-type CIs
were extracted from an extensive simulation of the limit distribution,
as no closed form expression is available.

\begin{figure}

\includegraphics{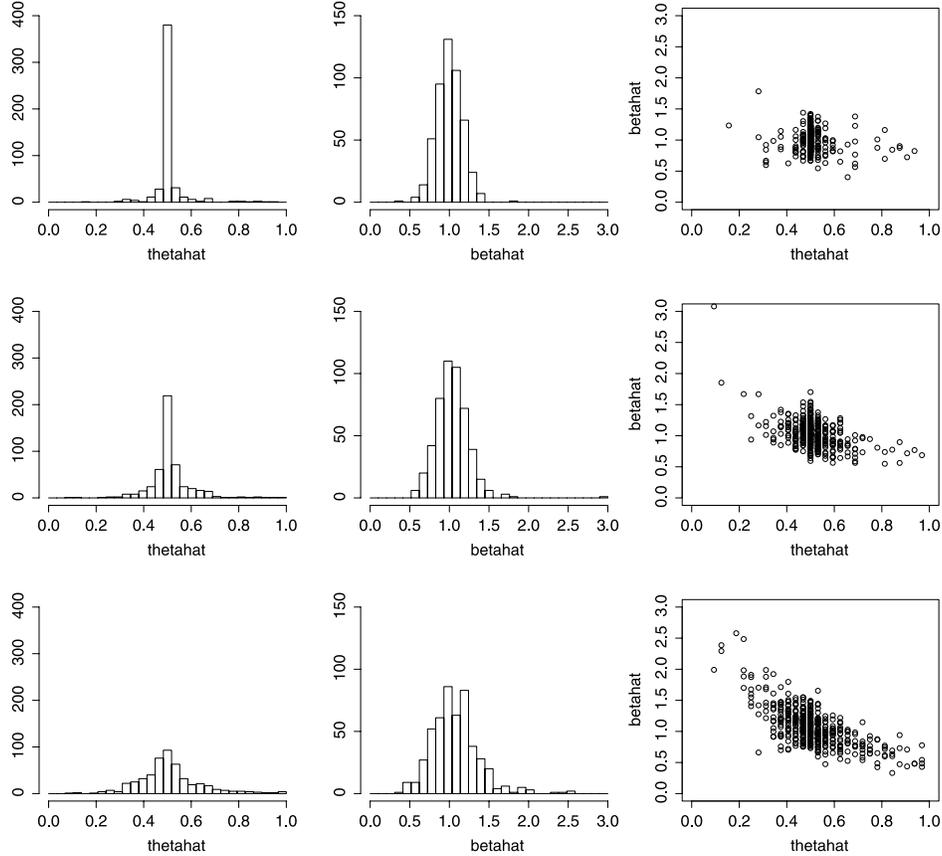}

\caption{Histograms and scatterplots of $\hat\theta_n$ and
$\hat\beta_n$ in the correctly specified case for $H=0.3$ (top row),
$H=0.5$ (middle row) and $H=0.7$ (bottom row), based on 500 samples of
size $n=20$.} \label{plot12}
\end{figure}

Table \ref{table:coverage_beta} reports the estimated coverage
probabilities and average lengths of nominal 95\% confidence intervals
for $\theta_0$ calculated using 500 independent samples. The data were
generated from the model (\ref{workmodel}), for $\alpha_0 = 0$,
$\beta_0 = 1$, $\theta_0 = 1/2$, $\varepsilon\sim N(0,\sigma^2)$ where
$\sigma= 0.3$ and $0.5$, the Hurst exponent $H = 0.3, 0.5, 0.7$ and
sample sizes $n=20$ and $40$. To calculate the least squares estimators
(\ref{lsq}), we restricted $\theta$ to a uniform grid of 101 points
in $[0,1]$; the fBm trajectories were generated over the same grid. The
fBm simulations were carried out in \textsf{R}, using the function \texttt{fbmSim}
from the \texttt{fArma} package, and via the Cholesky method of
decomposing the covariance matrix of $X$. Histograms and scatterplots
of $\hat\theta_n$ and $\hat\beta_n$ for $H=0.3, 0.5, 0.7$ when
$\sigma= 0.5$ are displayed in Figure \ref{plot12}.

%
\begin{table}
\caption{Monte Carlo results for coverage probabilities and average
widths of
nominal 95\% confidence intervals for $\theta_0$; data simulated from
the linear model with $\theta_0=1/2$,
$\alpha_0=0$ and $\beta_0=1$}\label{table:coverage_beta}
\begin{tabular*}{\tablewidth}{@{\extracolsep{\fill}}lcccccccc@{}}
\hline
& & &\multicolumn{2}{c}{\textbf{Wald-}$\bolds H$} &
\multicolumn{2}{c}{\textbf{R bootstrap}} &
\multicolumn{2}{c@{}}{\textbf{NP bootstrap}}\\[-4pt]
& & & \multicolumn{2}{c}{\hrulefill} &
\multicolumn{2}{c}{\hrulefill} &
\multicolumn{2}{c@{}}{\hrulefill}\\
$\bolds n$ & $\bolds\sigma$ & $\bolds H$ & \textbf{Cover}
& \textbf{Width} & \textbf{Cover} & \textbf{Width} & \textbf{Cover}
& \textbf{Width} \\
\hline
20 & 0.3 & 0.3 & 0.874 & 0.023 & 0.924 & 0.044 & 1.000 & 0.174 \\
& & 0.5 & 0.880 & 0.088 & 0.946 & 0.119 & 0.992 & 0.220 \\
& & 0.7 & 0.822 & 0.170 & 0.912 & 0.249 & 0.970 & 0.360 \\
& 0.5 & 0.3 & 0.806 & 0.129 & 0.912 & 0.211 & 0.998 & 0.410 \\
& & 0.5 & 0.852 & 0.256 & 0.924 & 0.333 & 0.988 & 0.487 \\
& & 0.7 & 0.834 & 0.352 & 0.938 & 0.510 & 0.962 & 0.591 \\
40 & 0.3 & 0.3 & 0.984 & 0.007 & 0.986 & 0.002 & 1.000 & 0.022 \\
& & 0.5 & 0.892 & 0.048 & 0.942 & 0.053 & 0.992 & 0.087 \\
& & 0.7 & 0.898 & 0.108 & 0.930 & 0.138 & 0.976 & 0.182 \\
& 0.5 & 0.3 & 0.900 & 0.039 & 0.928 & 0.054 & 0.998 & 0.149 \\
& & 0.5 & 0.908 & 0.134 & 0.950 & 0.165 & 0.990 & 0.251 \\
& & 0.7 & 0.856 & 0.229 & 0.946 & 0.332 & 0.962 & 0.386 \\
\hline
\end{tabular*}
\end{table}

%

In practice, $X$ can only be observed at discrete time points, so
restricting to a grid is the natural formulation for this example.
Indeed, the resolution of the observation times in the neighborhood of
$\theta_0$ is a limiting factor for the accuracy of $\hat\theta_n$,
so the grid resolution needs to be fine enough for the statistical
behavior of $\hat\theta_n$ to be apparent. For large sample sizes, a
very fine grid would be needed in the case of a small Hurst exponent
(cf. Theorem \ref{fbmcase}). 
Indeed, the histogram of $\hat\theta_n$ in the case $H=0.3$ (the
first plot in Figure \ref{plot12}) shows that the resolution of the
grid is almost too coarse to see the statistical variation, as the bin
centered on $\theta_0=1/2$ contains almost 80\% of the estimates. This
phenomenon is also observed in Table~\ref{table:coverage_beta} when $n
= 40$ and $\sigma= H = 0.3$. The average length of the CIs is smaller
than the resolution of the grid and, thus, we observe an over-coverage.
The two histograms of $\hat\theta_n$ for $H=0.5$ and $H=0.7$,
however, show increasing dispersion and become closer to bell-shaped as
$H$ increases.

Recall that, for simplicity, we pretend as if we know $H$, which should
be an advantage, yet the residual bootstrap has better performance
based on the results in Table \ref{table:coverage_beta}. We see that
usually the Wald-type CIs have coverage less than the nominal 95\%,
whereas the inconsistent nonparametric bootstrap method over-covers
with observed coverage probability close to 1. Accordingly, the average
lengths of the Wald-type CIs are the smallest, whereas those obtained
from the nonparametric bootstrap method are the widest. The behavior of
CIs obtained from the nonparametric bootstrap method also illustrates
the inconsistency of this procedure. A~similar phenomenon was also
observed in \cite{LM06} in connection with estimators that converge at
$n^{1/3}$-rate.

Despite the asymptotic independence of $\hat\theta_n$ and $\hat\beta
_n$, considerable correlation is apparent in the scatterplots in Figure
\ref{plot12}, with increasing negative correlation as $H$ increases;
note, however, that when $H=1$ there is a lack of identifiability of
$\theta$ and $\beta$, so the trend in the correlation as $H$
approaches 1 is to be expected in small samples.

%
\begin{figure}

\includegraphics{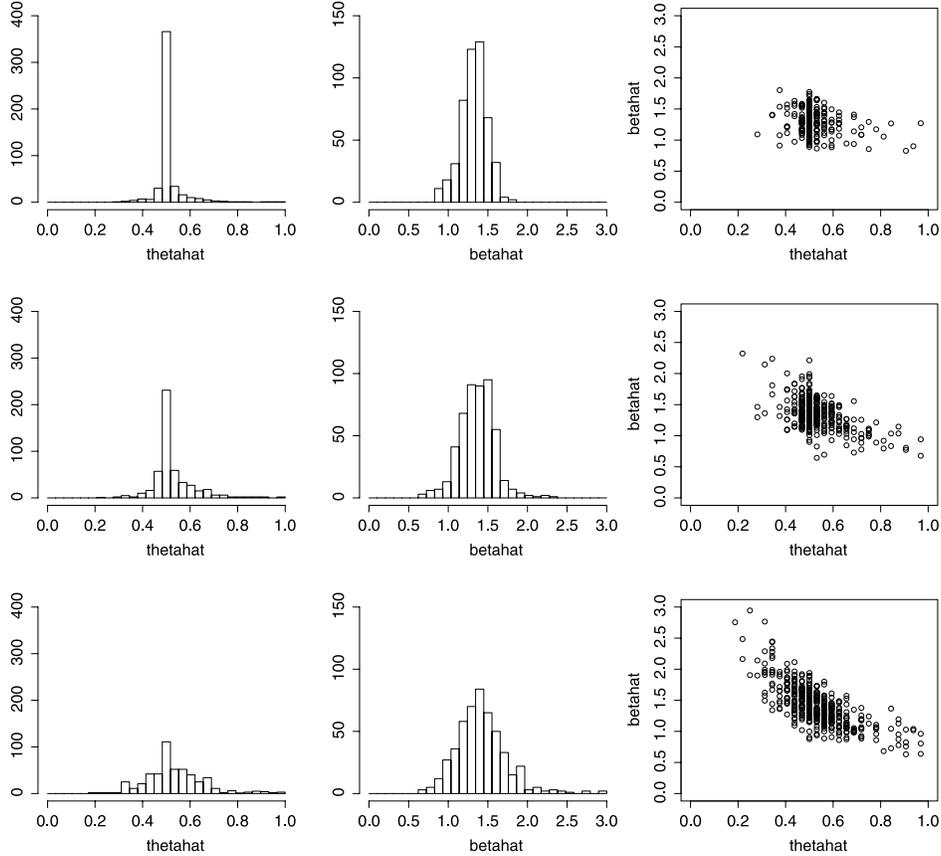}

\caption{Same as Figure \protect\ref{plot12} except in the partially
misspecified case.}\label{plot13}
\end{figure}

Next we consider a partially misspecified example, in which the data
are now generated from (\ref{partlinear}) by setting $f(t)=1/2$ and
$\theta=\theta_*=1/2$, but the other ingredients are unchanged from
the correctly specified example. The plots in Figure \ref{plot12}
correspond to those in Figure \ref{plot13}. The effect of
misspecification on $\hat\theta_n$ is a slight increase in dispersion
but no change in mean; the effect on $\hat\beta_n$ is a substantial
shift in mean along with a slight increase in dispersion.

\subsection{Gene expression example}\label{sec61}
Next we consider the gene
expression data mentioned in connection with Figure \ref{fig1},
to see how the residual bootstrap performs with such trajectories.
The trajectories consist of log gene expression levels from the breast
tissue of $n=40$ breast cancer patients, along a sequence of 518 loci
from chromosome 17. The complete gene expression data set is described
in Richardson et al. \cite{rich}. Although a continuous response is
not available for this particular data set, it is available in numerous
other studies
of this type; see the references mentioned in the \hyperref[sec1]{Introduction}.

%
\begin{figure}

\includegraphics{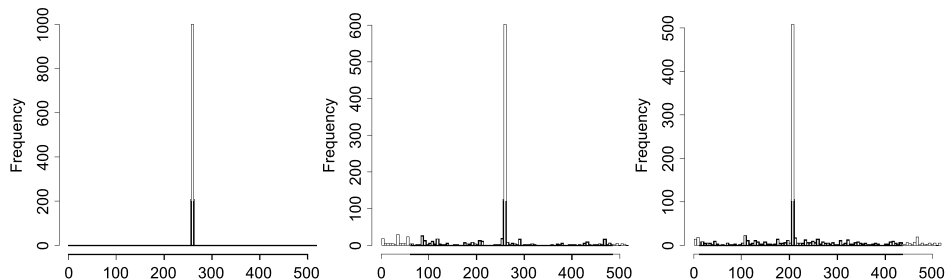}

\caption{Gene expression example: histograms of $\hat\theta_n^*$ based
on 1000 residual bootstrap samples and simulated responses with $\sigma
= 0.01$ (left), $\sigma= 0.03$ (middle) and $\sigma= 0.1$
(right).}\label{plotboot}
\end{figure}

To construct a scalar response, we generated $Y_i$ using the point
impact model (\ref{workmodel}) with $\alpha_0=0$ and $\beta_0=1$,
$\theta_0=0.5$ (corresponding to the\vspace*{1pt} position of 259 base pairs along
the chromosome)
and $\varepsilon\sim N(0,\sigma^2)$ for various values of $\sigma$.
As previously noted, the trajectories are very rough in this example
(with $H$ estimated to be about 0.1), which implies
a rapid rate of convergence for $\hat\theta_n$. We find that an
abrupt transition in the behavior of the residual bootstrap occurs
as $\sigma$ increases: for small $\sigma$, the residual bootstrap
estimates become degenerate at $\theta_0$ due to the
relatively coarse resolution; for moderately large $\sigma$, although
a considerable portion of the estimates are concentrated
at $\theta_0$, they become spread out over the 518 loci; for very
large $\sigma$, the estimates are more or less uniformly scattered
along the chromosome. Indeed, this is consistent with the behavior of
the Wald-type CI (\ref{eq:ci}) having width
proportional to $\sigma^{1/H}$, which blows up dramatically as $\sigma
$ increases when $H$ is small.

In Figure \ref{plotboot} we plot the bootstrap distribution of $\hat
\theta_n$ (obtained
from 1000 residual bootstrap samples in each case) for $\sigma=$ 0.01,
0.03 and 0.1. When $\sigma=0.01$, the bootstrap distribution is
degenerate at $\theta_0$; the resolution of the grid is
too course to see any statistical fluctuation in this case. When
$\sigma$ is moderate, namely, 0.03, although
the bootstrap distribution has a peak at $\theta_0$, the mass is widely
scattered and the resulting CI now covers almost
the entire chromosome. Further increasing the noise level causes the
bootstrap distribution to become even more dispersed and its mode moves
away from $\theta_0$; the sample size of 40 is
now too small for the method to locate the neighborhood of
$\theta_0$.

\section{Concluding remarks}\label{Discuss}

In this paper we have developed a point impact functional linear
regression model for use with trajectories as predictors of a
continuous scalar response. It is expected that the proposed approach
will be useful when there are sensitive time points at which the
trajectory has an effect on the response. We have derived the rates of
convergence and the explicit limiting distributions of the least
squares estimator of such a parameter in terms of the Hurst exponent for
fBm trajectories. We also established the validity of the residual
bootstrap method for obtaining CIs for sensitive time points, avoiding
the need to estimate the Hurst exponent. In addition, we have developed
some results in the misspecified case in which the data are generated
partially or completely from a standard functional linear model, and in
the two-sample setting.

Although for simplicity of presentation we have assumed that the
trajectories are fBm, it is clear from the proofs that
it is only local properties of the trajectories in the neighborhood of
the sensitive time point that drive the theory, and thus
the validity of the confidence intervals. The consistency of the least
squares estimator is of course needed,
but this could be established under much weaker assumptions (namely,
uniform convergence of
the empirical criterion function and the existence of a well-separated
minimum of the limiting
criterion function; cf. \cite{vw}, page 287).

Exploiting the fractal behavior of the trajectories plays a crucial
role in developing confidence intervals based on
the least squares estimator of the sensitive time point, in contrast to
standard functional
linear regression where it is customary to smooth the predictor
trajectories prior to fitting the model (\cite{rs}, Chapter 15).
Our approach does not require any preprocessing of the trajectories
involving a choice of smoothing
parameters, nor any estimation of nuisance parameters (namely, the
Hurst exponent).
On the other hand, functional linear regression is designed with
prediction in mind, rather than
interpretability, so in a sense the two approaches are complimentary.
The tendency
of functional linear regression to over-smooth a point impact (see
\cite{lmck} for detailed discussion)
is also due to the use of a roughness penalty on the regression
function; the smoothing parameter is usually
chosen by cross-validation, a criterion that optimizes for predictive
performance.

Our model naturally extends to allow multiple sensitive time points,
and any model selection procedure having the oracle property
(such as the adaptive lasso) could be used to estimate the number of
those sensitive time points.
The bootstrap procedure for the (unregularized) least squares estimator
can then be adapted to provide individual CIs around each time point,
although developing theoretical justification would be challenging.
Other challenging problems would be to develop bootstrap procedures
that are suitable for the two-sample
problem and for the misspecified model settings.

It should be feasible to carry through much of our program for certain
types of diffusion processes driven by fBm, and also for processes
having jumps. In the case of piecewise constant trajectories that have
a single jump, the theory specializes to an existing type of
change-point analysis \cite{kqs}. Other possibilities include L\'{e}vy
processes (which have stationary independent increments) and
multi-parameter fBm. It should also be possible to develop versions of
our results in the setting of censored survival data (e.g., Cox
regression). Lindquist and McKeague \cite{lmck} recently studied point
impact generalized linear regression models in the case that $X$ is
standard Brownian motion and we expect that our approach can be
extended to such models as well.

\section{Proofs}\label{proofs}
To avoid measurability problems and for simplicity of notation, we will
always use outer expectation/probability, and denote them by $E$ and $P$;
$E^*$ and $P^*$ will denote bootstrap conditional
expectation/probability given the data.

We begin with the proof of Theorem \ref{fbmcase}.
The strategy is to establish (a) consistency, (b) the rate of
convergence, (c) the weak convergence of a suitably localized version
of the criterion function, and (d) apply the argmax (or argmin)
continuous mapping theorem.

\subsection{Consistency}\label{sec81}
We start with some notation. Let $m_\eta(Y,X) \equiv[Y - \alpha-
\beta X(\theta)]^2$ and let $\mathbb{M}_n(\eta) \equiv\mathbb{P}_n
m_\eta= \frac{1}{n} \sum_{i=1}^n [Y_i - \alpha- \beta X_i(\theta
)]^2$, where $\mathbb{P}_n$ denotes the expectation with respect to
the empirical measure of the data. 
Let
%
%
\begin{eqnarray}\label{eq:M_eta}\qquad
\mathbb{M}(\eta) & \equiv& P m_\eta= (\alpha_0 - \alpha)^2 + P[\{
\beta_0 X(\theta_0) - \beta X(\theta)\}^2] + \sigma^2 \nonumber\\
& = & (\alpha_0 - \alpha)^2 + \sigma^2 + (\beta_0 - \beta)^2 P
[X^2(\theta_0)] + \beta^2 P [X(\theta_0) - X(\theta)]^2 \\
& &{} +2 \beta(\beta_0 - \beta) P[X(\theta_0)\{X(\theta_0) -
X(\theta)\}]. \nonumber
\end{eqnarray}
First observe that $\mathbb{M}(\eta)$ has a unique minimizer at $\eta
_0$ as $P[\beta X(\theta) \ne\beta_0 X(\theta_0)] > 0$, for all
$(\beta,\theta) \in\mathbb{R} \times(0,1)$ with $(\beta,\theta)
\ne(\beta_0, \theta_0)$.

Using the fBm covariance formula (\ref{covar}),
%
%
\begin{eqnarray}\label{eq:P_M}
\mathbb{M}(\eta) - \mathbb{M}(\eta_0) & = & (\alpha_0 - \alpha)^2
+ (\beta_0 -\beta)^2
|\theta_0|^{2H} + \beta^2 |\theta_0 - \theta|^{2H} \nonumber\\
& &{} + \beta(\beta_0 - \beta) \{|\theta_0|^{2H} + |\theta_0 -
\theta|^{2H} - |\theta|^{2H}\} \nonumber\\[-8pt]\\[-8pt]
& = & (\alpha_0 - \alpha)^2 + (\beta_0 -\beta)^2 |\theta_0|^{2H} +
\beta\beta_0 |\theta_0 - \theta|^{2H} \nonumber\\
& &{} + \beta(\beta_0 - \beta) \{|\theta_0|^{2H} - |\theta|^{2H}\}.
\nonumber
\end{eqnarray}

To show that $\hat\eta_n$ is a consistent estimator of $\eta_0$,
first note that $\hat\eta_n$ is uniformly tight. Also
notice that $\mathbb{M}(\eta)$ is continuous and has a unique minimum
at $\eta_0$, and, thus, by Theorem 3.2.3(i) of \cite{vw}, it is
enough to show that $\mathbb{M}_n \stackrel{P} {\rightarrow} \mathbb
{M} $
uniformly
on each compact subset $K$ of
$\Xi= \mathbb{R}^2 \times[0,1]$, and that $\mathbb{M}$
has a well-separated minimum in the sense that $\mathbb{M}(\eta_0) <
\inf_{\eta\notin G} \mathbb{M}(\eta)$ for every open set $G$ that
contains $\eta_0$. That $\mathbb{M}$ has a well-separated minimum can
be seen from the form of the expression in (\ref{eq:P_M}). For the
uniform convergence, we need to show that the class $\mathcal{F} = \{
m_\eta\dvtx\eta\in K\}$ is $P$-Glivenko Cantelli ($P$-GC). Using GC
preservation properties (see Corollary 9.27 of \cite{kbook}), it is
enough to show that $\mathcal{G} = \{B_H(h) \equiv X(\theta_0 + h) -
X(\theta_0) \dvtx h \in[-1,1]\}$ is $P$-GC. Note that almost all
trajectories of $X$ are Lipschitz of any order strictly less than $H$,
in the sense of (\ref{lip})
in Lemma \ref{lem2.1} below. Thus, the bracketing number $ N_{[\cdot
]}(\varepsilon, \mathcal{G}, L_1(Q)) < \infty$ and $\mathcal{G}$ is
$P$-GC, by Theorems 2.7.11 and 2.4.1 of \cite{vw}.

\subsection{Rate of convergence}\label{sec82}
We will apply a result of van der Vaart and Wellner (\cite{vw},
Theorem 3.2.5) to obtain a lower bound on the rate of convergence of the
M-estimator $\hat\eta_n$. Setting $\tilde d(\eta,\eta_0) = \max\{
|\alpha- \alpha_0|, |\beta- \beta_0|, |\theta- \theta_0|^{H}\}$,
the first step is to show that
%
%
\begin{equation}\label{eq:M_eta-M_eta0}
\mathbb{M}(\eta) - \mathbb{M}(\eta_0) \gtrsim\tilde d^2(\eta,\eta_0)
\end{equation}
in a neighborhood of $\eta_0$. Here $\gtrsim$ means that the right-hand side
is bounded above by a (positive) constant times the left-hand side.
Note that $|\theta_0|^{2H} - |\theta|^{2H}$ has a bounded derivative
in $\theta\in[\delta,1]$, where $\delta> 0$ is fixed, so for such
$\theta$ we have
%
%
\begin{eqnarray}\label{eq:SimM_eta}
&&\beta(\beta- \beta_0) \{|\theta_0|^{2H} - |\theta|^{2H}\} \nonumber\\
&&\qquad \ge
- |\beta| |\beta_0 - \beta| C |\theta_0 - \theta| \nonumber\\[-8pt]\\[-8pt]
&&\qquad = - [|\beta| C |\theta_0 - \theta|^{1-H} ] |\beta_0 - \beta|
|\theta_0 - \theta|^{H} \nonumber\\
&&\qquad \ge - c(\theta,\beta) \{ (\beta_0 - \beta)^2 + |\theta_0 -
\theta|^{2H} \},\nonumber
\end{eqnarray}
where $C$ is the bound on the derivative, $c(\theta,\beta) = |\beta|
C |\theta_0 - \theta|^{1-H}/2$, and we used the inequality $|ab| \le
(a^2 + b^2)/2$. As $\beta_0 \ne0$ and $0 < \theta_0 < 1$, by
combining (\ref{eq:P_M}) and (\ref{eq:SimM_eta}), suitably grouping
terms, and noting that $c(\theta,\beta)$ can be made arbitrarily
small by restricting $\eta$ to a sufficiently small neighborhood of
$\eta_0$, there exist $c_1 > 0$ and $c_2 > 0$ such that
\[
\mathbb{M}(\eta) - \mathbb{M}(\eta_0) \ge(\alpha_0 - \alpha)^2 +
c_1 (\beta_0 -\beta)^2 + c_2 |\theta_0 - \theta|^{2H},
\]
which shows that (\ref{eq:M_eta-M_eta0}) holds.

Let $\mathcal{M}_{\delta} \equiv\{m_{\eta} - m_{\eta_0}\dvtx\tilde
d(\eta,\eta_0) < \delta\}$, where $\delta\in(0,1]$. Note that
%
%
\begin{eqnarray}\label{eq:m_eta}\qquad
m_\eta- m_{\eta_0} & = & (\alpha^2 - \alpha_0^2) + \beta
^2[X^2(\theta) - X^2(\theta_0)] + (\beta^2 - \beta_0^2) X^2(\theta
_0) \nonumber\\
&&{} - 2Y(\alpha- \alpha_0) - 2 \beta Y [X(\theta) - X(\theta_0)] -
2 (\beta- \beta_0) Y X(\theta_0) \nonumber\\[-8pt]\\[-8pt]
&&{} + 2 \alpha\beta[X(\theta) - X(\theta_0)] + 2 \alpha X(\theta
_0) (\beta- \beta_0) \nonumber\\
&&{} + 2 \beta_0 X(\theta_0) (\alpha- \alpha_0).\nonumber
\end{eqnarray}
This shows that $\mathcal{M}_{\delta}$ has envelope
%
%
\begin{eqnarray}\label{eq:envelop}\hspace*{5pt}
M_\delta(Y,X)& \equiv& {(2 |\alpha_0| + \delta) \delta+ (|\beta
_0| + \delta)^2 \sup_{|\theta- \theta_0|^H < \delta} }|X^2(\theta)
- X^2(\theta_0)| \nonumber\\
&&{} + X^2(\theta_0) \delta(2 |\beta_0| + \delta) + 2 |Y| \delta
\nonumber\\
&&{} + 2 |Y| {(|\beta_0| + \delta) \sup_{|\theta- \theta_0|^H <
\delta}} |X(\theta) - X(\theta_0)| \\
&&{} + 2 |X(\theta_0)||Y| \delta+ 2 (|\alpha_0| + \delta) {(|\beta_0| +
\delta) \sup_{|\theta- \theta_0|^H < \delta}} |X(\theta) -
X(\theta_0)| \nonumber\\
&&{} + 2 (|\alpha_0| + \delta) |X(\theta_0)| \delta+ 2 |\beta_0|
|X(\theta_0)| \delta.\nonumber
\end{eqnarray}
Using a maximal inequality for fBm (Theorem 1.1 of \cite{nov}), we have
%
%
\begin{equation}
\label{self-sim}
E \Bigl[ \sup_{|\theta- \theta_0|^H < \delta} |X(\theta) - X(\theta
_0)|^q \Bigr] \lesssim\delta^{q}
\end{equation}
for any $q>0$. Then, using (A3) in conjunction with H\"{o}lder's inequality
(cf. the proof of Lemma \ref{lem2.1}),
all nine terms in (\ref{eq:envelop}) can be shown to have second moments
bounded by $\delta^2$ (up to a constant) and, thus, $E M_\delta^2
\lesssim\delta^2$.

The following lemma shows that $m_\eta$ is ``Lipschitz in parameter''
and, consequently, that the bracketing entropy integral $J_{[\cdot
]}(1,\mathcal{M}_\delta,L^2(P))$ is uniformly bounded as a function of
$\delta\in(0,1]$; see \cite{vw}, page 294. Without loss of
generality, to simplify notation, we assume that $\alpha= 0$ and
$\beta= 1$, and state the lemma with $\theta$ as the only parameter.
\begin{lemma}\label{lem2.1}
If \textup{(A1)} and \textup{(A3)} hold and $0<\alpha<H$, there is a random variable $L$
with finite second moment such that
%
%
\begin{equation}\label{lip}
|m_{\theta_1} -m_{\theta_2}| \le L |\theta_1 - \theta_2|^{\alpha}
\end{equation}
for all $\theta_1, \theta_2 \in[0,1]$ almost surely.
\end{lemma}
\begin{pf} The trajectories of fBm are Lipschitz of any order $\alpha
<H$ in the sense that
%
%
\begin{equation}\label{lipbm}
|X(t) -X(s)| \le\xi|t-s|^\alpha\qquad \forall t,s\in[0,1]
\end{equation}
almost surely, where $\xi$ has moments of all orders; this is a
consequence of the proof of Kolmogorov's continuity theorem; see
Theorem 2.2 of Revuz and Yor~\cite{ry}. Noting that $ m_\theta
(X,Y)=(Y - X(\theta))^2$, we then have
\[
|m_{\theta_1} -m_{\theta_2}|\le C|X(\theta_1) - X(\theta_2)| \le
L|\theta_1-\theta_2|^\alpha,
\]
where $C=2({\sup_\theta}|X(\theta)|+|Y|)$ and $L=C\xi$.
Here $L$ has a finite second moment:
\[
E L^2\le\{EC^{2p}\}^{1/p}\{ E\xi^{2 q}\}^{1/q}<\infty
\]
by H\"{o}lder's inequality for $1/p + 1/q = 1$ with
$p=1+\delta/2$ and $\delta>0$ comes from the moment condition (A3).
\end{pf}

Using a maximal inequality from \cite{vw} (see page 291), we then have
\[
E_P \| \mathbb{G}_n \|_{\mathcal{M}_\delta}
\lesssim J_{[\cdot]}(1, \mathcal{M}_\delta, L_2(P)) (E M_\delta^2)^{1/2}
\lesssim\delta
\]
for all $\delta\in(0,1]$, where $\mathbb{G}_n=\sqrt n(\mathbb{P}_n-P)$,
and it follows that $\tilde d(\hat\eta_n,\eta_0) = O_P(1/\sqrt n)$
by Theorem 3.2.5 of \cite{vw}.

\subsection{Localizing the criterion function}\label{sec83}
To simplify\vspace*{1pt} notation, let $r_n^{-1} \mathbf{h} \equiv (h_1/\sqrt
{n}, h_2/\sqrt{n}, n^{-1/(2H)} h_3)$, for $\mathbf{h} = (h_1,h_2,h_3)
\in\mathbb{R}^3$. Then
%
%
\begin{equation}\label{basicdecomp}
\zeta_n = \mathop{\arg\min}_{\mathbf{h}} [\mathbb{M}_n(\eta_0 + r_n^{-1}
\mathbf{h})-
\mathbb{M}_n(\eta_0)]
\end{equation}
and we can write the expression in the square brackets after
multiplication by $n$ as the sum of an empirical process and a drift term:
%
%
\begin{equation}\label{decom}
\mathbb{G}_n \bigl[\sqrt{n}(m_{\eta_0 + r_n^{-1}
\mathbf{h}} - m_{\eta_0}) \bigr]
+ n[\mathbb{M}(\eta_0 + r_n^{-1} \mathbf{h}) - \mathbb{M}(\eta_0)].
\end{equation}
First consider the empirical process term, and note that
\begin{eqnarray*}
m_{\eta_0 + r_n^{-1} \mathbf{h}} & = & \bigl[Y - (\alpha_0 + n^{-1/2}
h_1) - (\beta_0+ n^{-1/2} h_2) X\bigl(\theta_0 + n^{-1/(2H)} h_3\bigr) \bigr]^2 \\
& = & \biggl[\varepsilon- \biggl\{ \frac{h_1}{\sqrt{n}} + \biggl( \beta_0 + \frac
{h_2}{\sqrt{n}} \biggr) X\bigl(\theta_0 + n^{-1/(2H)} h_3\bigr) - \beta_0 X(\theta
_0) \biggr\} \biggr]^2,
\end{eqnarray*}
so we obtain
%
%
\begin{eqnarray}\label{eq:n(Pn-P)}\quad
\sqrt{n} [m_{\eta_0 + r_n^{-1} \mathbf{h}} - m_{\eta_0} ] &=& \sqrt
{n} \biggl[ \frac{h_1}{\sqrt{n}} + \biggl( \beta_0 + \frac{h_2}{\sqrt{n}} \biggr)
\frac{\mathbb{B}(h_3)}{\sqrt{n}} + \frac{h_2}{\sqrt{n}} X(\theta
_0) \biggr]^2 \nonumber\\[-8pt]\\[-8pt]
&&{}- 2 \varepsilon\biggl[ h_1 + \biggl(\beta_0 + \frac{h_2}{\sqrt{n}} \biggr) \mathbb
{B}(h_3) + h_2 X(\theta_0) \biggr],\nonumber
\end{eqnarray}
where $\mathbb{B}(h_3) \equiv\sqrt{n} [X(\theta_0 + n^{-1/(2H)} h_3)
- X(\theta_0)] \stackrel{d}{=} B_H(h_3)$ (as a process in $h_3$).

The result of applying $\mathbb{G}_n$ to the first term on the right-hand
side of the above display
gives a term of order $o_P(1)$ uniformly in $\mathbf{h} \in[-K,K]^3$,
for each
$K>0$.
This is seen
by applying the maximal inequality from \cite{vw}, page 291, as used
above; here the
class of functions $\mathcal{F}_n$ in question
is bounded by the envelope function
\[
F_n = 3 \sqrt{n} \biggl\{ {\frac{K^2}{n} + \biggl( \beta_0 + \frac{K}{\sqrt{n}}
\biggr)^2 \sup_{|h_3| \le K}} \frac{\mathbb{B}^2(h_3)}{n} + \frac{K^2}{n}
X^2(\theta_0) \biggr\},
\]
for which $P F_n^2 = o(1)$ and $J_{[\cdot]} (1,\mathcal
{F}_n,L_2(P))<\infty$; cf. the proof of Lemma \ref{lem2.1}.
Hence, we just need to consider the second term.
To determine the limit distribution of the empirical process term in
(\ref{decom}), it thus suffices to show that
%
%
\begin{equation}
\label{tosee}
\mathbb{G}_n[(\varepsilon, \varepsilon\mathbb{B}(h_3), \varepsilon
X(\theta
_0))] \stackrel{d}{ \to} (\sigma Z_1,\sigma B_H(h_3),\sigma Z_2)
\end{equation}
in $\mathbb{R}\times C[-K,K] \times\mathbb{R}$, where $Z_1,Z_2$ are
i.i.d. $N(0,1)$ and independent of the fBm~$B_H$.
For the second component above, notice that since $\varepsilon$ is
independent of~$\mathbb{B}$,
%
%
\begin{equation}
\label{simple}
\mathbb{G}_n [\varepsilon\mathbb{B}(h_3)]\stackrel{d}{=} B_H(h_3)
\Biggl({1\over
n}\sum_{i=1}^n\varepsilon_i^2 \Biggr)^{1/2}
\stackrel{d}{ \to} \sigma B_H(h_3)
\end{equation}
in $C[-K,K]$.
The asymptotic independence of the three components of (\ref{tosee})
is a consequence of
\begin{eqnarray*}
\operatorname{Cov}(\varepsilon, \varepsilon\mathbb{B}(h_3)) & = & \sigma^2
E[\mathbb{B}(h_3)] = 0, \\
\operatorname{Cov}(\varepsilon, \varepsilon X(\theta_0)) & = & \sigma^2
E[X(\theta
_0)] = 0, \\
\operatorname{Cov}(\varepsilon\mathbb{B}(h_3), \varepsilon X(\theta_0))
& = &
\sigma^2 \biggl[\frac{\sqrt{n}}{2} ( |\theta_0 + n^{-1/2H} h_3|^{2H} -
|\theta_0|^{2H} ) - \frac{h_3}{2\sqrt{n}} \biggr],
\end{eqnarray*}
which tends to zero uniformly in $h_3 \in[-K,K]$, using the assumption
$H < 1$.

It just remains to find the limit of the drift term in (\ref{decom}).
Using (\ref{eq:P_M}), it is given by
\begin{eqnarray*}
&& h_1^2 + h_2^2 |\theta_0|^{2H} + ( \beta_0 + n^{-1/2} h_2 ) \beta
_0 |h_3|^{2H}\\
&&\quad{} + h_2 ( \beta_0 + n^{-1/2} h_2 ) \bigl[ \sqrt{n} \{ |\theta_0|^{2H} -
| \theta_0 + n^{-1/2H} h_3 |^{2H} \} \bigr] \\
&&\qquad \rightarrow h_1^2 + h_2^2 |\theta_0|^{2H} + \beta_0^2 |h_3|^{2H}
\end{eqnarray*}
uniformly in $\mathbf{h} \in[-K,K]^3$. Combining this with the limit
distribution of the first term in (\ref{decom}), we get from (\ref
{basicdecomp}) and the argmin continuous mapping theorem that
\begin{eqnarray*}
\zeta_n & \stackrel{d}{\rightarrow} & \mathop{\arg\min}_{\mathbf{h}} \bigl[ -2
\sigma\bigl(Z_1 h_1 + \beta_0 B_H(h_3) + h_2 |\theta_0|^H Z_2\bigr) \\
&&\hspace*{69.6pt}{} + (h_1^2 + h_2^2 |\theta_0|^{2H} + \beta_0^2 |h_3|^{2H}) \bigr] \\
& \stackrel{d}{=} & \biggl[ \sigma Z_1, |\theta_0|^{-H} \sigma Z_2, \mathop
{\arg
\min}_{h_3} \biggl\{ 2 \frac{\sigma}{|\beta_0|} B_H(h_3) + |h_3|^{2H} \biggr\} \biggr].
\end{eqnarray*}
This completes the proof of Theorem \ref{fbmcase}.

\subsection{\texorpdfstring{Proof of Theorem \protect\ref{thm:InConsBoots}}{Proof of Theorem 3.1}}\label{sec84}

We prove the result by the method of contradiction. Before giving the proof, we
state a general lemma that can be useful in studying bootstrap
validity. The lemma can be proved easily using characteristic
functions; see also Sethuraman \cite{sethu}
and Theorem 2.2 of Kosorok \cite{k08}.
\begin{lemma}\label{lemma:independence} Let $W_n$ and $W_n^*$ be
random vectors in $\mathbb{R}^l$ and $\mathbb{R}^k$, respectively; let
$Q$ and $Q^*$ be distributions on the Borel sets of $\mathbb{R}^l$ and
$\mathbb{R}^k$, and let $\mathcal{F}_n$ be $\sigma$-fields for which
$W_n$ is $\mathcal{F}_n$-measurable. If $W_n$ converges in
distribution to $Q$ and the conditional distribution of $W_n^*$ given
$\mathcal{F}_n$ converges (in distribution) in probability to $Q^*$,
then $(W_n,W_n^*)$ converges in distribution to $Q \times Q^*$.
\end{lemma}

The basic idea of the proof of the theorem now is to assume that
$\Delta_n^* \stackrel{d}{\rightarrow} \Delta^*$ in probability,
where $\Delta^*$ has the same distribution as $\Delta$. Therefore,
$\Delta_n^* \stackrel{d}{\rightarrow} \Delta^*$ unconditionally
also. We already know that $\Delta_n \stackrel{d}{\rightarrow}
\Delta$ from Theorem \ref{fbmcase}. By Lemma \ref
{lemma:independence} applied with $W_n = \Delta_n$, $W_n^* = \Delta
_n^*$ and ${\mathcal F}_n = \sigma((Y_1,X_1),(Y_2,X_2),\break\ldots,
(Y_n,X_n))$, we can show that $(\Delta_n,\Delta_n^*)$ converges
unconditionally to a product measure and, thus, $\Delta_n + \Delta
_n^* \stackrel{d}{\rightarrow} \Delta+ \Delta^*$. Thus, $n^{1/(2H)}
(\hat\theta_n^* - \theta_0) \equiv\Delta_n + \Delta_n^*$
converges unconditionally to a tight limiting distribution which has
twice the variance of $\Delta$.

Using arguments along the lines of those used in the proof of Theorem
\ref{fbmcase}, we can show that
\[
n^{1/(2H)} (\hat\theta_n^* - \theta_0) \stackrel{d}{\rightarrow}
\mathop{\arg\min}_{t \in\mathbb{R}} \bigl\{ 2 \sigma\beta_0 \bigl(B_H(t) + B_H^*(t)\bigr)
+ \beta_0^2 |t|^{2H} \bigr\} \equiv\Delta^{**},
\]
where $B_H^*$ is another independent fBm with Hurst exponent $H$. Using
properties of fBm, we see that
\[
\Delta^{**} \stackrel{d}{=} \biggl( \sqrt{2} \frac{\sigma}{|\beta_0|}
\biggr)^{1/H} \mathop{\arg\min}_{t \in\mathbb{R}} \{B_H(t) + |t|^{2H}/2 \}
\stackrel{d}{=} 2^{1/(2H)} \Delta.
\]
Thus, the variance of the limiting distribution of $n^{1/(2H)} (\hat
\theta_n^* - \theta_0)$ is $2^{1/H} > 2$ times the variance of
$\Delta$, which is a contradiction.

\subsection{\texorpdfstring{Proof of Theorem \protect\ref{ConsBoots}}{Proof of Theorem 3.2}}\label{sec85}

The bootstrap sample is $\{(Y_i^*,X_i), i=1,\ldots,n\}$, where the
$Y_i^*$ are defined in (\ref{eq:DefYi*}). Letting $\mathbb
{M}^*_n(\eta) \equiv{\mathbb{P}^*_n} m_\eta= \frac{1}{n} \sum
_{i=1}^n [Y_i^* - \alpha- \beta X_i(\theta)]^2$, the bootstrap
estimates are
%
%
\begin{equation}\label{eq:Cov}
\hat\eta_n^* = (\hat\alpha_n^*, \hat\beta_n^*, \hat\theta_n^*)
\equiv\mathop{\arg\min}_{\eta\in\Xi} \mathbb{M}_n^*(\eta).
\end{equation}
%
We omit the rate of convergence part of the proof,
and concentrate on establishing the limit distribution. Also, to keep
the argument simple, we will
assume that $\hat\eta_n \to\eta_0$ a.s., but a subsequence argument
can be used to bypass this assumption. Note that
%
%
\begin{equation}\label{Decomp}\qquad
\zeta_n^* = \mathop{\arg\min}_{\mathbf{h} \in\mathbb{R}^3} \{ n(\mathbb
{P}_n^* - P_n) [m_{\hat\eta_n + r_n^{-1} \mathbf{h}} - m_{\hat\eta
_n}] + nP_n [m_{\hat\eta_n + r_n^{-1} \mathbf{h}} - m_{\hat\eta
_n}]\},
\end{equation}
where $P_n$ is the probability measure generating the bootstrap sample.
Consider the first term within the curly brackets. Using a similar
calculation as in (\ref{eq:n(Pn-P)}),
%
%
\begin{equation}\label{eq:BootsTerm1}
\sqrt{n}(m_{\hat\eta_n + r_n^{-1} \mathbf{h}} - m_{\hat\eta_n}) =
- 2 \varepsilon^*[h_1 + \hat\beta_n \hat{\mathbb{B}}(\hat\theta_n,
h_3) + h_2 X(\hat\theta_n)] + A_n,
\end{equation}
where $\hat{\mathbb{B}}(\theta, t) \equiv\sqrt{n} [X(\theta+
n^{-1/(2H)} t) - X(\theta)]$, with the dependence on $n$ suppressed
for notational convenience, and $a_n\equiv\sqrt n(\mathbb{P}_n^*-P_n)
A_n =
o_P(1)$ uniformly in $\mathbf{h} \in[-K,K]^3$. Then, using (\ref
{eq:BootsTerm1}),
%
%
\begin{eqnarray}\label{eq:BootsTerm2}
&& \sqrt{n}(\mathbb{P}_n^* - P_n)\bigl[\sqrt{n}(m_{\hat\eta_n +
r_n^{-1} \mathbf{h}} - m_{\hat\eta_n})\bigr] \nonumber\\
&&\qquad = - \sqrt{n}(\mathbb{P}_n^* - P_n) [\varepsilon^*\{h_1 + \hat
\beta
_n \hat{\mathbb{B}}(\hat\theta_n, h_3) + h_2 X(\hat\theta_n)\}] +
a_n\\
&&\qquad \stackrel{d}{\rightarrow} -2 \sigma\bigl(Z_1 h_1 + \beta_0 B_H(h_3)
+ h_2 |\theta_0|^H Z_2\bigr)\nonumber
\end{eqnarray}
in $C[-K,K]$, a.s., where $Z_1,Z_2$ are i.i.d. $N(0,1)$ that are
independent of $B_H$.

To prove (\ref{eq:BootsTerm2}), first note that $P_n [\varepsilon^*\{h_1
+ \hat\beta_n \hat{\mathbb{B}}(\hat\theta_n, h_3) + h_2 X(\hat
\theta_n)\}] = 0$, as the $X_i$ are fixed and the $\varepsilon_i^*$ have
mean zero under $P_n$. We will need the following properties of $\hat
{\mathbb{B}}(\hat\theta_n, t)$, proved at the end:
%
%
\begin{eqnarray}\label{eq:PropB}
\frac{1}{n} \sum_{i=1}^n \hat{\mathbb{B}}_i(\hat\theta_n, t)
&\stackrel{P}{\rightarrow}& 0,\qquad \frac{1}{n} \sum_{i=1}^n \hat{\mathbb
{B}}_i(\hat\theta_n, t) X_i(\hat\theta_n)\stackrel{P}{\rightarrow
} 0, \nonumber\\[-8pt]\\[-8pt]
\frac{1}{n} \sum_{i=1}^n \hat{\mathbb{B}}_i(\hat\theta_n, s)
\hat{\mathbb{B}}_i(\hat\theta_n, t) &\stackrel{P}{\rightarrow}&
C_H(s,t),\nonumber
\end{eqnarray}
uniformly for $|s|,|t| <K$, where $C_H(s,t) $ is the covariance
function (\ref{covar}) of fBm. Now considering (\ref{eq:BootsTerm2}),
by simple application of the Lindeberg--Feller theorem,
it follows that
\[
\sqrt{n} \mathbb{P}_n^* [\varepsilon^* h_1] \stackrel
{d}{\rightarrow}
h_1 N(0, \sigma^2),\qquad \sqrt{n} \mathbb{P}_n^* [\varepsilon^* h_2
X(\hat
\theta_n)] \stackrel{d}{\rightarrow} h_2 N(0, |\theta_0|^{2 H}
\sigma^2),
\]
a.s. in $C[-K,K]$. Next consider
$ \sqrt{n} \mathbb{P}_n^* [\varepsilon^* \hat{\mathbb{B}}(\hat
\theta
_n, t)] $. The finite-dimensional convergence and tightness of this
process follow from Theorems 1.5.4 and 1.5.7 in \cite{vw} using the
properties of $\hat{\mathbb{B}}(\hat\theta_n, t)$ stated in (\ref
{eq:PropB}). The asymptotic independence of the terms under
consideration also follows using (\ref{eq:PropB}) via a similar
calculation as in (\ref{eq:Cov}).

To study the drift term in (\ref{Decomp}), note that
%
%
\begin{eqnarray}
P_n m_\eta& = & \frac{1}{n} \sum_{i=1}^n P_n[Y_i^* - \alpha- \beta
X_i(\theta)]^2 \nonumber\\
& = & \frac{1}{n} \sum_{i=1}^n \frac{1}{n} \sum_{j=1}^n [\hat
\alpha_n + \hat\beta_n X_i(\hat\theta_n) + (\hat\varepsilon_j -
\bar{\varepsilon}_n)- \alpha- \beta X_i(\theta)]^2
\nonumber\\[-8pt]\\[-8pt]
& = & \frac{1}{n} \sum_{i=1}^n [ (\hat\alpha_n - \alpha) + (\hat
\beta_n - \beta) X_i(\hat\theta_n) + \beta\{ X_i(\hat\theta_n) -
X_i(\theta) \} ]^2 \nonumber\\
&&{} + \frac{1}{n} \sum_{j=1}^n (\hat\varepsilon_j - \bar
{\varepsilon
}_n)^2. \nonumber
\end{eqnarray}
Simple algebra then simplifies the drift term to
%
%
\begin{eqnarray}\label{eq:BootsTerm3}
& & \sum_{i=1}^n \biggl\{ \frac{h_1}{\sqrt{n}} + \frac{h_2}{\sqrt{n}}
X_i(\hat\theta_n) + \frac{\hat{\mathbb{B}}_i(\hat\theta_n,
h_3)}{\sqrt{n}} \biggl( \hat\beta_n + \frac{h_2}{\sqrt{n}} \biggr) \biggr\}^2
\nonumber\\
&&\qquad = h_1^2 + \frac{h_2^2}{n} \sum_{i=1}^n X_i^2(\hat\theta_n) + \biggl(
\hat\beta_n + \frac{h_2}{\sqrt{n}} \biggr)^2 \frac{1}{n} \sum_{i=1}^n
\hat{\mathbb{B}}_i(\hat\theta_n, h_3)^2 \nonumber\\
&&\qquad\quad{} + 2 \frac{h_1 h_2}{n} \sum_{i=1}^n X_i(\hat\theta_n) + 2 h_1 \biggl(
\hat\beta_n + \frac{h_2}{\sqrt{n}} \biggr) \frac{1}{n} \sum_{i=1}^n
\hat{\mathbb{B}}_i(\hat\theta_n, h_3) \\
&&\qquad\quad{} + 2 h_2 \biggl( \hat\beta_n + \frac{h_2}{\sqrt{n}} \biggr) \frac{1}{n}
\sum_{i=1}^n \hat{\mathbb{B}}_i(\hat\theta_n, h_3) X_i(\hat\theta
_n) \nonumber\\
&&\qquad \stackrel{P}{\rightarrow} h_1^2 + h_2^2 |\theta_0|^{2H} + \beta
_0^2 |h_3|^{2H}\nonumber
\end{eqnarray}
uniformly on $[-K,K]$, where we have used the properties of $\hat
{\mathbb{B}}(\hat\theta_n, h_3)$ in (\ref{eq:PropB}) and
\[
\Biggl| \frac{1}{n} \sum_{i=1}^n X_i(\hat\theta_n) \Biggr|
\le{\sup_{\theta}}
|(\mathbb{P}_n - P) X(\theta)| \stackrel{P}{\rightarrow} 0.
\]
Thus, combining (\ref{Decomp}), (\ref{eq:BootsTerm2}) and (\ref
{eq:BootsTerm3}), we get $ \zeta_n^* \stackrel{d}{\rightarrow} \zeta$
in probability.

It remains to prove (\ref{eq:PropB}). We only prove the last part, the
other parts being similar. For fixed $K >0 $, consider the function
class
\[
\mathcal{F}_n = \{\hat{\mathbb{B}}(\theta,s) \hat{\mathbb
{B}}(\theta,t)\dvtx\theta\in[0,1], |s| < K, |t| < K \},
\]
which has a uniformly bounded bracketing entropy integral, and envelope
\[
F_n = {\sup_{\theta, |s| < K, |t| < K } }|\hat{\mathbb{B}}(\theta,s)
\hat{\mathbb{B}}(\theta,t)|\le n^{\alpha'/H} K^{2(H-\alpha')} \xi^2
\]
from the Lipschitz property (\ref{lipbm}) of order $\alpha=H-\alpha
'$, where $0< \alpha' < H/2 $ and $\xi$ has finite moments of all
orders. Then
\begin{eqnarray*}
&& P \Biggl\{ \sup_{|s|, |t| < K} \Biggl|\frac{1}{n} \sum_{i=1}^n \hat{\mathbb
{B}}_i(\hat\theta_n, t) \hat{\mathbb{B}}_i(\hat\theta_n, s) -
C_H(s,t) \Biggr| > \varepsilon\Biggr\} \\
&&\qquad \le P \Bigl\{ {\sup_{f \in\mathcal{F}_n} }|(\mathbb{P}_n - P) f | >
\varepsilon\Bigr\} \le\frac{1}{\varepsilon} E \Bigl[ {\sup_{f
\in\mathcal{F}_n}}
|(\mathbb{P}_n - P) f | \Bigr] \\
&&\qquad \lesssim \frac{1}{\varepsilon\sqrt{n}} J_{[\cdot]}(1,\mathcal
{F}_n,L_2(P)) (E F_n^2)^{1/2} \lesssim n^{\alpha'/H - 1/2} \rightarrow0,
\end{eqnarray*}
where we use a maximal inequality in Theorem 2.14.2 of \cite{vw}.

\subsection*{Remark} The failure of the nonparametric bootstrap
can be explained from the behavior of the drift term in (\ref{Decomp}).
In the nonparametric bootstrap, we need to find the conditional limit
of $n\mathbb{P}_n[m_{\hat\eta_n + r_n^{-1} \mathbf{h}} - m_{\hat
\eta_n}]$ given the data, but observe
that $\sqrt n \mathbb{P}_n$ applied to the second term of (\ref{eq:n(Pn-P)})
fails to converge
in probability. However, when bootstrapping residuals, the drift term
in (\ref{Decomp})
becomes $n{P}_n[m_{\hat\eta_n + r_n \mathbf{h}} - m_{\hat\eta_n}]$,
and $\sqrt n P_n$ applied to the second term in (\ref{eq:n(Pn-P)})
vanishes, so the drift term now converges in
probability, as seen in (\ref{eq:BootsTerm3}).

\subsection{\texorpdfstring{Proof of Theorem \protect\ref{miss}}{Proof of Theorem 4.1}}\label{sec86}

The consistency of $\hat\theta_n$ follows using a
Glivenko--Cantelli argument for the function class
${\mathcal F}\equiv\{m_\theta(X,Y)=[Y-X(\theta)]^2\dvtx\theta\in
[0,1]\}$ and the existence of a well-separated
minimum for $\mathbb{M}$; cf. the proof of Theorem~\ref{fbmcase}.
Note that $\theta_0$ is the unique solution of the normal equation
$\mathbb{M}'(\theta)=0$
and $\mathbb{M}''(\theta_0)>0$, so
%
%
\begin{equation}
\label{metric}
\mathbb{M}(\theta)-\mathbb{M}(\theta_0)\gtrsim d^2(\theta,\theta_0)
\end{equation}
for all $\theta$ in a neighborhood of $\theta_0$, where $d$ is the
usual Euclidean distance.
The envelope function $M_\delta= {\sup_{|\theta-\theta_0|<\delta
}}|m_\theta-m_{\theta_0}|$ for $\mathcal{M}_\delta\equiv
\{m_\theta-m_{\theta_0}\dvtx\theta\in[0,1]\}$ has $L^2$-norm of order
$ \delta^H$, from (\ref{self-sim}), so
Theorem 3.2.5 of \cite{vw} applied with $\phi_n(\delta) =\delta^H$
gives rate $r_n=n^{1/(4-2H)}$
with respect to Euclidean distance.

Now write
$\hat h_n\equiv r_n (\hat\theta_n -\theta_0)=\arg\min_{h\in
\mathbb{R}
}\tilde{\mathbb{M}}_n(h)$,
where
%
%
\begin{equation}
\label{local}
\tilde{\mathbb{M}}_n(h)=r_n^2[\mathbb{M}_n(\theta
_0+h/r_n)-\mathbb{M}_n(\theta_0)],\qquad
h\in\mathbb{R}.
\end{equation}
This gives
%
%
\begin{equation}\label{decompFBM}\qquad
\tilde{\mathbb{M}}_n(h)=n^{-H/(4-2H)}\mathbb{G}_n[Z_n(h)^2] -
2\mathbb{G}_n [W Z_n(h)]
+ \tfrac{1}{2}\mathbb{M}''(\theta_0)h^2+A_n,
\end{equation}
where $A_n=o(1)$ uniformly in $h\in[-K,K]$, for any $K>0$, and
\begin{eqnarray*}
W&\equiv&\int_0^1 f(t)X(t) \,dt -X(\theta_0)+\varepsilon,\\
Z_n(h) &\equiv&
n^{H/(4-2H)}[X(\theta_0 + h/r_n) - X(\theta_0)].
\end{eqnarray*}
Note that $Z_n(h)=_d B_H(h)$ as processes, so, by Donsker's theorem,
the first term in (\ref{decompFBM}) converges to zero in probability
uniformly over $[-K,K]$. For the second term, we claim that
%
%
\begin{equation}\label{keystep}
\mathbb{G}_n [W Z_n(h) ]\stackrel{d}{\rightarrow} a B_H(h)
\end{equation}
as processes in $C[-K,K]$, where $a^2= E(W^2)$. Application
of the argmin continuous mapping theorem will then complete the proof.

To prove (\ref{keystep}), for simplicity, we just give the detailed
argument in the Brownian motion case, with
$B=B_{1/2}$ denoting two-sided Brownian motion.
Consider the decomposition
%
%
\begin{equation}\label{decompG}
\mathbb{G}_n [W Z_n(h) ] = \mathbb{G}_n [(W-W_\eta)Z_n(h) ] +\mathbb
{G}_n [W_\eta Z_n(h) ],
\end{equation}
where
%
%
\begin{equation}\quad
W_{\eta}=\int_{\theta_0-\eta}^{\theta_0+\eta}f(t)X(t) \,dt
+\bigl(X(\theta_0+\eta)-X(\theta_0)\bigr)\bigl(
F(1)-F(\theta_0+\eta)\bigr),
\end{equation}
$F(\theta)=\int_0^\theta f(t) \,dt$, and $\eta>0$ is sufficiently
small so that
$|\theta_0\pm\eta|<1$. Splitting the range of integration
for the first term in $W$ into three intervals, and using the
integration by parts
formula (for semimartingales)
over the intervals $[0,\theta_0-\eta]$ and $[\theta_0+\eta,1]$, we get
\begin{eqnarray*}
W-W_\eta&=&\int_0^{\theta_0-\eta} \bigl(F(\theta_0-\eta)-F(t)\bigr)\, dX(t)+
\int_{\theta_0+\eta}^1 \bigl(F(1)-F(\theta_0+\eta)\bigr) \,dX(t)
\\
&&{} + \varepsilon+X(\theta_0)\bigl(F(1)-F(\theta_0+\eta)-1\bigr),
\end{eqnarray*}
which implies, by the independent increments property, that $W-W_\eta$
is independent of $Z_n(h)$ for $|h|<\eta n^{1/3}$.
Using the same argument
as in proving (\ref{simple}), it follows
that
\[
\mathbb{G}_n [(W-W_\eta) Z_n(h) ]\stackrel{d}{\rightarrow} a_\eta B(h)
\]
as processes in $C[-K,K]$, for each fixed $\eta>0$, where
\[
a_\eta^2= E(W-W_\eta)^2\to E(W^2)=E \biggl[\int_0^1 f(t)X(t) \,dt -X(\theta
_0) \biggr]^2+ \sigma^2\equiv a^2
\]
as $\eta\to0$. Clearly,
$ a_\eta B(h)\stackrel{d}{\rightarrow} aB(h)$ in $C[-K,K]$ as $\eta\to0$.
If we show that the last term in (\ref{decompG}) is asymptotically negligible
in the sense that, for every $M>0$ and $\delta>0$,
%
%
\begin{equation}
\label{bill}
\lim_{\eta\to0}\limsup_{n\to\infty}P \Bigl({\sup_{|h|<M}}|\mathbb
{G}_n[W_\eta
Z_n(h)]|>\delta\Bigr)=0,
\end{equation}
this will complete the proof in view of Theorem 4.2 in \cite{bill}.
Theorem 2.14.2 in \cite{vw} gives
\[
E \Bigl[{\sup_{|h|<M}}|\mathbb{G}_n[W_\eta Z_n(h)]| \Bigr]\lesssim J_{[\cdot
]}(1,\mathcal
{F},L^2(P)) (E F^2)^{1/2},
\]
where $J_{[\cdot]}(1,\mathcal{F},L^2(P))$ is the bracketing entropy
integral of the class of functions $\mathcal{F}=\mathcal{F}_{n,\eta}
=\{ W_\eta Z_n(h)\dvtx|h|<M\}$, and
$F=F_{n,\eta}$ is an envelope function for $\mathcal{F}$. We can take
$F=|W_\eta| {\sup_{|h|<M}} |Z_n(h)|$. By the Cauchy--Schwarz inequality,
\[
E(F^2)\le( E W_\eta^4 )^{1/2} \Bigl( {E \sup_{|h|<M} }|B(h)|^4
\Bigr)^{1/2}\lesssim\eta M,
\]
where we have used (\ref{self-sim}). The bracketing entropy integral
can be shown to be uniformly bounded (over all $\eta>0$ and $n$)
using the Lipschitz property (\ref{lipbm}). The previous two displays
and Markov's inequality then lead to
\[
\limsup_{n\to\infty}P \Bigl({\sup_{|h|<M}}|\mathbb{G}_n[W_\eta
Z_n(h)]|>\delta
\Bigr)\lesssim\sqrt{\eta M}/\delta,
\]
which implies (\ref{bill}) and establishes (\ref{keystep}).


To establish (\ref{keystep}) for general fBm, we apply Theorem 2.11.23
of \cite{vw} to the class of measurable functions $\mathcal{F}_n =\{
f_{n,h}\dvtx|h|<M\}$, where $f_{n,h}(X,\varepsilon)= WZ_n(h)$ and
$M>0$ is
fixed. Direct computation using the covariance of fBm shows that the
sequence of covariance functions of $f_{n,h}$ converges pointwise to
the covariance function of $aB_H(h)$, and the various other conditions
can be shown to be satisfied using similar arguments to what we have
seen already.

\subsection{\texorpdfstring{Proof of (\protect\ref{2samp})}{Proof of (13)}}\label{sec87}

The key step involving the localization of the criterion function again
relies on the self-similarity of fBm $B_H$:
\begin{eqnarray*}
n_1^{1/(2H)} (\hat\theta_n-\theta_0)&=&\mathop{\arg\max}_{h} (\mathbb
{P}_n^1-\mathbb{P}
_n^2)\bigl[X\bigl(\theta_0+n_1^{-1/(2H)}h\bigr)-X(\theta_0)\bigr]\\
&\stackrel{d}=& \mathop{\arg\max}_{h}\bigl\{ \bigl(\mathbb{G}_n^1-\sqrt\rho\mathbb{G}
_n^2\bigr)[B_H(h)] \\
& &\hspace*{37.9pt}{} + n_1 \bigl(\mathbb{M}\bigl(\theta_0+n_1^{-1/(2H)}h\bigr)-\mathbb{M}(\theta
_0)\bigr)\bigr\}\\
& \stackrel{d}{\rightarrow}& \mathop{\arg\max}_{h}\bigl\{ \bigl(1+\sqrt\rho\bigr) B_H(h)
-c|h|^{2H}\bigr\},
\end{eqnarray*}
where $\mathbb{G}_n^j=\sqrt{n_j}(\mathbb{P}_n^j-P_j)$ is the
empirical process for
the $j$th sample.

\section*{Acknowledgments}
The authors thank Moulinath Banerjee and Davar Khoshnevisan for helpful
discussions.

\printaddresses

\end{document}